\documentclass{article}
\usepackage{latexsym,amssymb}

\title{Bounding the order of complex linear groups and permutation groups with selected composition factors}

\author{Geoffrey R. Robinson}

\begin{document}

\maketitle

\begin{abstract}
Originally motivated by questions of P. Etingof (related to growth rates of tensor powers in symmetric tensor categories, see [4]), we obtain (as usual, modulo Abelian normal subgroups) general bounds of exponential nature on the order of finite subgroups $G$ of ${\rm GL}(n,\mathbb{C})$ with prescribed properties.  Our most general result of this nature is :

\medskip
\noindent {\bf Theorem:} \emph{Let ${\bf P}$ be a property of (isomorphism types of) finite groups which is inherited both by subgroups, and by homomorphic images. Then one and only one of the following is true:}

\medskip
\noindent i) \emph{ Every finite group has property ${\bf P}.$}

\medskip
\noindent ii)  \emph{ There is a real number $c({\bf P})$  such that whenever $n$ is a positive integer and $G$ is a finite subgroup of ${\rm GL}(n,\mathbb{C})$ which has property ${\bf P}$, then $G$ has an Abelian normal subgroup $A$ with $[G:A] \leq c({\bf P})^{n-1}.$}

\medskip
We prove this theorem as a consequence of a more precise result on properties ${\bf Q}$ which are inherited by normal subgroups and homomorphic images, where (using the Classification of the Finite Simple Groups), we are able to determine exactly when such a constant $c({\bf Q})$ exists, and to give explicit upper bounds for it when it does. Note that since the symmetric group ${\rm Sym}_{n+1}$ is isomorphic to a subgroup of ${\rm GL}(n,\mathbb{C})$ for all $n,$ the two options in the Theorem are indeed mutually exclusive.

\medskip
The two questions from P. Etingof concerned groups of order prime to $p$, and groups with Abelian Sylow $p$-subgroups, respectively, where $p$ is a specified prime. We are able to provide the optimal upper bound in almost every case.  By placing the questions in this general framework, we are able to indicate numerous other directions to apply our methods.

\end{abstract}

\section{Introduction, definitions, and statement of main results}

\medskip
We recall that a \emph{section} of the finite group $G$ is a group of the form $H/K,$ where $H$ is isomorphic to a subgroup of $G$ and $K \lhd H.$  A \emph{normal section} is a section $H/K$ with $H \lhd G.$ By an \emph{overgroup} of a finite group $G$, we mean a group $M$ such that $G$ is isomorphic to a subgroup of $M$.

\medskip
Throughout, we let ${\bf P}$ be a property of finite groups which is inherited by all sections (hereditary for brevity) , and we let ${\bf Q}$ be a property which is inherited by all normal sections (normally hereditary for brevity).  We assume throughout that having property ${\bf P}$ or ${\bf Q}$ is invariant under isomorphism.

\medskip
Note that ${\bf Q}$ being normally hereditary is equivalent to requiring that ${\bf Q}$ is inherited both by normal subgroups and by homomorphic images. The chosen property ${\bf P}$ is also normally hereditary, but the chosen property ${\bf Q}$ need not be hereditary in general.

\medskip
We will say that the normally hereditary property ${\bf Q}$ is \emph{extant} if at least one finite group has property  ${\bf Q}$, and we will say that the extant property ${\bf Q}$ is \emph{irredundant} if there is at least one  finite group which does not have property ${\bf Q}.$ 

\medskip
As usual, when $G$ is a finite group, we let $F(G)$ denote the Fitting subgroup of $G$,  that is, the subgroup of $G$ generated by all nilpotent normal subgroups of $G$ (which is also the unique largest nilpotent normal subgroup of $G$).

\medskip
When ${\bf Q}$ is an extant normally hereditary property of finite groups, we define the sets of real numbers $A({\bf Q})$, $B({\bf Q})$ and 
$C({\bf Q})$ as follows:

\medskip
$A({\bf Q}) = \{ [G:F(G)]^{\frac{1}{n-1}} : n \geq 2, G$ \emph{is a finite subgroup of} ${\rm GL}(n,\mathbb{C})$ \emph{with property}  ${\bf Q}  \}.$

\medskip
$B({\bf Q}) = \{ |G|^{\frac{1}{n-1}} : n \geq 2, G \leq {\rm Sym}_{n},$ \emph{G has property} ${\bf Q} \}.$

\medskip
$C({\bf Q}) = \{ [G:A]^{\frac{1}{n-1}} : n \geq 2, G$ \emph{is a finite subgroup of} ${\rm GL}(n,\mathbb{C})$ \emph{with property} ${\bf Q}$,  \emph{ $A$ is an Abelian normal subgroup of $G$ of maximal order} $\}.$

\medskip
Our first general result, which makes use of the classification of the finite simple groups, is:

\medskip
\noindent {\bf Theorem 1.1:} \emph{ Let ${\bf P}$ be an extant irredundant hereditary property of finite groups. Then the sets $A({\bf P})$, $B({\bf P})$, and $C({\bf P})$ are all bounded above.}

\medskip
\emph{ Hence letting $\alpha({\bf P})$, $\beta({\bf P})$ and 
and $\gamma({\bf P})$ denote the least upper bounds of $A({\bf P}),B({\bf P})$ and $C({\bf P})$ respectively, we have :}

\medskip
\noindent i) \emph{ Whenever $n$ is a positive integer and $G$ is a finite subgroup with property ${\bf P}$ of ${\rm GL}(n,\mathbb{C}),$
then we have $[G:F(G)] \leq \alpha({\bf P})^{n-1}$. Furthermore, $\alpha({\bf P})$ is the smallest real number with this property.}

\medskip
\noindent ii) \emph{ Whenever $n$ is a positive integer and $G$ is a finite subgroup with property ${\bf P}$ of the symmetric group ${\rm Sym}_{n},$
then we have $|G| \leq \beta({\bf P})^{n-1}.$ Furthermore, $\beta({\bf P})$ is the smallest real number with this property.}

\medskip
\noindent iii) \emph{ Whenever $n$ is a positive integer and $G$ is a finite subgroup with property ${\bf P}$ of ${\rm GL}(n,\mathbb{C}),$
then $G$ has an Abelian normal subgroup $A$ with $[G:A] \leq \gamma({\bf P})^{n-1}.$ Furthermore, $\gamma({\bf P})$ is the smallest real number with this property.}

\medskip
Theorem 1.1 has no serious content when only finitely many groups have property ${\bf P},$ and it requires explicit upper bounds for the constants involved to find applications of other than theoretical interest. The applications to follow are derived from appropriate choices of the property ${\bf P}$, and we will obtain explicit bounds for  these constants in due course.

\medskip
Theorem 1.1 is a consequence of the following theorem on normally hereditary properties ${\bf Q}$. We will see in due course that
this bound can be improved for some properties ${\bf Q}.$  For expository convenience, we will henceforth always consider the alternating group 
${\rm Alt}_{2}$ to be the trivial group.

\medskip
\noindent {\bf Theorem 1.2:} \emph{ Let ${\bf Q}$ be an extant normally hereditary property of finite groups. Then the sets $A({\bf Q})$, $B({\bf Q})$, and $C({\bf Q})$ are all bounded above if and only if there are only finitely many alternating groups with property ${\bf Q}.$}

\medskip
\emph{Furthermore, when these sets are all bounded above, and we let $m$ denote the largest integer such that ${\rm Alt}_{m}$ has property ${\bf Q},$ then we
 have:}

\medskip
\noindent \emph{ i)  If $m \leq 151,$ then $${\rm max}\{ \alpha({\bf Q}), \beta({\bf Q}),\gamma({\bf Q}) \} \leq 60.$$  }

\medskip
\noindent \emph{ ii) If $m >151,$ then $$\left(\frac{m!}{2} \right)^{\frac{1}{m-1}} \leq {\rm min} \{\alpha({\bf Q}),\beta({\bf Q}),\gamma({\bf Q}) \}$$ and 
$${\rm max}\{ \alpha({\rm Q}), \beta({\bf Q}),\gamma({\bf Q}) \} \leq \left(m!\right)^{\frac{1}{m-2}}.$$ }

\medskip
That Theorem 1.2 implies Theorem 1.1 follows from the next simple Lemma, whose proof we can  give now.

\medskip
\noindent {\bf Lemma 1.3:} \emph{Let ${\bf P}$ be an extant irredundant hereditary property of finite groups. Then only finitely many alternating groups have property ${\bf P}.$}

\medskip
\noindent {\bf Proof:} We recall that for all $t \in \mathbb{N},$ the symmetric group ${\rm Sym}_{t}$ has a natural embedding into the alternating group ${\rm Alt}_{t+2}.$ Together with Cayley's theorem, this tells us that every finite group $G$ is isomorphic to a subgroup
of ${\rm Alt}_{|G|+2}.$

\medskip
By assumption, there is a finite group $H$ which does not have property ${\bf P}$. Since ${\bf P}$ is inherited by sections, and, in particular,
by subgroups, we see that no overgroup of $H$ has property ${\bf P}.$ Setting $k = |H|+2,$ we conclude that no alternating group ${\rm Alt}_{n}$ 
with $n \geq k$ has property ${\bf P},$ since $H$ is isomorphic to a subgroup of ${\rm Alt}_{n}$ for all $n \geq k.$ Hence only finitely many alternating groups have property ${\bf P}$.

\medskip
\noindent {\bf Remark 1.4 :} If we can determine the largest integer $m$ such ${\rm Alt}_{m}$ has property ${\bf P}$, then we may apply Theorem 1.2 to obtain explicit bounds for Theorem 1.1. The proof of Lemma 1.3 gives an upper bound for $m$ when we can locate a single finite group $H$ which fails to have property ${\bf P}$. More generally than in that proof, if the chosen group $H$ has a corefree subgroup $K$ of index $t$, then $H$ embeds in ${\rm Alt}_{t}$ when $H = O^{2}(H)$, and embeds in ${\rm Alt}_{t+2}$ in any case. Hence in all cases, ${\rm Alt}_{t+2}$ does not have property ${\bf P}.$ 

\medskip
The questions of P. Etingof give two situations when Theorem 1.1 may be applied to obtain explicit bounds. If we fix a prime $p$, and the property 
${\bf P}$ is the property of having order prime to $p$, then ${\bf P}$ is clearly hereditary, extant, and irredundant. For $p$ an odd prime, $p-1$ is the largest  integer $m$ such that ${\rm Alt}_{m}$ has property ${\bf P}$, while for $p = 2$, only the alternating groups ${\rm Alt}_{2}$ and ${\rm Alt}_{3}$ have property ${\bf P}.$

\medskip
If we instead take ${\bf P}$ to be the property of having Abelian Sylow $p$-subgroups (for our chosen prime $p$), then the largest integer $m$ such that the alternating group ${\rm Alt}_{m}$ has property ${\bf P}$ is  $p^{2}-1$ if $p$ is odd, and is $p^{2}+1$ if $p = 2.$ Hence Theorem 1.2 may be applied, and we obtain:

\medskip
\noindent {\bf Corollary 1.5:} \emph{Let $p$ be a prime :}

\medskip
\noindent \emph{ i) If $p < 157$, and $G$ is a finite subgroup of ${\rm GL}(n,\mathbb{C})$ of order prime to $p$,
then $G$ has an Abelian normal subgroup of index at most $60^{n-1}.$}

\medskip
\noindent \emph{ ii) If $p \geq  157,$ and $G$ is a finite subgroup of ${\rm GL}(n,\mathbb{C})$ of order prime to $p$,
then $G$ has an Abelian normal subgroup of index at most $(p-1)!^{\frac{n-1}{p-3}}.$ This upper bound is attained 
when $G = {\rm Sym}_{p-1}$ and $n = p-2.$}

\medskip
\noindent \emph{ iii) If $p  \leq 11$, and $G$ is a finite subgroup of ${\rm GL}(n,\mathbb{C})$ with Abelian Sylow $p$-subgroups,
then $G$ has an Abelian normal subgroup of index at most $60^{n-1}.$}

\medskip
\noindent \emph{ iv) If $p > 11,$ and $G$ is a finite subgroup of ${\rm GL}(n,\mathbb{C})$ with Abelian Sylow $p$-subgroups,
then $G$ has an Abelian normal subgroup of index at most $(p^{2}-1)!^{\frac{n-1}{p^{2}-3}}.$ This upper bound is attained when $G = {\rm Sym}_{p^{2}-1}$ and $n = p^{2}-2 .$}

\medskip
\noindent {\bf Remark 1.6:} The bound in parts i)-ii)  can be improved when $p = 2$, but this was implicitly dealt with in [10], as discussed in [4].
Because of the example $G = {\rm SL}(2,5)$ in its two-dimensional representation, the bound is optimal when $p > 5.$

\medskip
The bound in parts iii)-iv) can be improved via our methodology to $60^{\frac{n-1}{2}}$ when $p = 2$, since ${\rm SL}(2,5)$ does not have Abelian Sylow $2$-subgroups, but we omit the details. For odd primes $p$, the bound of parts iii)-iv) is optimal.

\medskip
If we fix a prime $p$ and a positive integer $k$, we may define the extant, hereditary, irredundant, property ${\bf P}$ to be that of having a Sylow $p$-subgroup of nilpotency class at most $k$. Then there are always finite groups which do not have property ${\bf P}$, and we may conclude from Theorems 1.1 and  1.2 that  $\alpha({\bf P}),\beta({\bf P})$ and $\gamma({\bf P})$ all exist. 

\medskip
For example, the nilpotency class of a Sylow $p$-subgroup of ${\rm Sym}_{n}$ is $p^{r}$ when $p^{r+1} \leq n < p^{r+2}$, so 
we may easily find the largest integer $m$ such that a Sylow $p$-subgroup  of ${\rm Alt}_{m}$ has nilpotency class at most $k$ for any integer $k,$ and then we may apply Theorem 1.2 to obtain explicit upper bounds for $\alpha({\bf P}),\beta({\bf P})$ and $\gamma({\bf P}).$

\medskip
Clearly, there are numerous other properties ${\bf P}$ and ${\bf Q}$ to which we can apply Theorems 1.1 and 1.2, but we single out one more for explicit mention, since (at least for $n > 150$) it refines a theorem of M.J. Collins [2], which in turn improved a theorem of B. Weisfeiler whose proof remained incomplete and unpublished due to his unfortunate disappearance (see however [12],[13]). 

\medskip
We indicate here how to derive its proof, given Theorem 1.2  (together with an observation from Section 3 below, that the sequence $ (m!^{\frac{1}{m-2}})$ increases for $m \geq 5).$

\medskip
\noindent {\bf Corollary 1.7:} \emph{Let $H$ be a finite subgroup of ${\rm GL}(n,\mathbb{C})$. If $H$ has any non-Abelian alternating composition factor, let $m \geq 5$ be the largest integer such that ${\rm Alt}_{m}$ is a composition factor of $H$. Then either $H$ has an Abelian normal subgroup of index at most $60^{n-1},$ or else $n+1 \geq  m > 151,$ and $H$ has an Abelian normal subgroup $A$ with
$$\frac{m!}{2} \leq [H:A] \leq \left( m! \right)^{\frac{n-1}{m-2}} \leq (n+1)!.$$}

\medskip
\noindent {\bf Proof:} We define the property ${\bf Q}$ by stipulating that the finite group $G$ has property ${\bf Q}$ if and only if either $G$ is trivial, or else every composition factor of $G$ occurs as a composition factor of $H$. The property ${\bf Q}$ is normally hereditary (but not hereditary, in general), extant, and irredundant. Furthermore, since $H$ has only finitely many composition factors, it is certainly the case that only finitely many alternating groups have property ${\bf Q},$ so that Theorem 1.2 may be applied. Let $A$ be an Abelian normal subgroup of $H$ of minimal index.

\medskip
If $H$ has no non-Abelian alternating composition factor, then we are done by Theorem 1.2. Otherwise, let ${\rm Alt}_{m}$ be the largest non-Abelian alternating composition factor of $H$. If $m \leq 151,$ then we are done by Theorem 1.2. 

\medskip
Suppose then that $m > 151.$ Then ${\rm Alt}_{m}$ is a composition factor of $H/A$, so that $[H:A] \geq \frac{m!}{2}.$   On the other hand, Theorem 1.2 implies that $$[H:A] \leq (m!)^{\frac{n-1}{m-2}}.$$  Hence $m! \leq 2(m!)^{\frac{n-1}{m-2}}$ and so
$(m!)^{m-2} \leq  2^{m-2}(m!)^{n-1}.$ If $m > n+1$, this yields $m! \leq  2^{m-2},$ a contradiction. Hence $m \leq n+1$
so that $$(m!)^{\frac{n-1}{m-2}} \leq (n+1)!^{\frac{n-1}{n+1-2}} = (n+1)!.$$

\medskip
\section{Further description of results and proofs} 

\noindent
We will presently see that for many properties ${\bf P}$ and ${\bf Q}$, we can improve the constant $60$ in Theorem 1.2, which leads (for example) to improvements in Corollary 1.5 for some primes $p$, and to Corollary 1.7 for many groups $H$. Since Theorem 1.1 is a consequence of Theorem 1.2, we concentrate our efforts on proving (sometimes strengthened versions of) Theorem 1.2.

\medskip
The methods here are similar in principle to the paper [10], though we are working in much more generality here. In [10], (translating to current terminology), we were working with normally hereditary properties ${\bf Q}$  such that the finite group $G$ had property ${\bf Q}$ if and only if all prime divisors of $|G|$ are greater than or equal to $p$ for some specified odd prime $p$. We proved that the least upper bounds $\alpha({\bf Q}),\beta({\bf Q}),$ and  $\gamma({\bf Q})$ all exist, and we found estimates for them which decreased as $p$ increased. The discussion in Section 1.1 of [10] is relevant to, and may shed some light on, our work here. Our primary focus, both here and in [10], is on finite subgroups $G$ of ${\rm GL}(n,\mathbb{C}).$ However, while it is advantageous for purposes of induction and self-containment to use a unified methodology to treat permutation groups and linear groups together, the bounds we obtain here in passing for permutation groups almost never  improve existing bounds, whereas they sometimes did in [10] for solvable permutation groups of odd order. 

\medskip
For example, A. Maroti has already proved in [9] that if $G$ is a primitive subgroup of ${\rm Sym}_{n}$ which does not contain ${\rm Alt}_{n},$
then $|G| < 3^{n},$ and $|G| < 2^{n}$ if $n >24$. For primitive subgroups of ${\rm Sym}_{n}$, there are often 
polynomial-type bounds if alternating composition factors are avoided, and something similar is true for finite primitive subgroups $G$ of ${\rm GL}(n,\mathbb{C})$ (bounding $[G:F(G)]$ or $[G:A]$ when $A$ is an Abelian normal subgroup of $G$ of minimal index). For example, in his preliminary draft [13] and article [12], Weisfeiler used bounds of the  form $n^{a}$ or $n^{b \log{n}}$ for constants $a$ and $b$ ( when alternating composition factors are excluded). However, such sub-exponential  bounds are inevitably destroyed (both for permutation groups and linear groups) by natural wreath product constructions, even in solvable and nilpotent groups, so do not usually persist in imprimitive groups.

\medskip
We recall that a group $H$ is said to be \emph{almost simple} if it has a unique minimal normal subgroup which is non-Abelian and simple. A \emph{component} of a group $G$ is a quasi-simple subnormal subgroup $L$ of $G$ (that is to say, $L \lhd \lhd G$,  $L = L^{\prime}$ and $L/Z(L)$ is simple). Distinct components of $G$ centralize each other, and centralize the Fitting subgroup $F(G)$. The central product of all components of $G$ is denoted by $E(G)$, and the \emph{generalized Fitting subgroup} $F^{\ast}(G)$ is defined to be the central product $F(G)E(G).$ We always have $C_{G}(F^{\ast}(G)) = Z(F(G))$ and it is always the case that $G/F^{\ast}(G)$ is isomorphic to a subgroup of ${\rm Out}_{G}(F^{\ast}(G))$, the subgroup of the outer automorphism of $F^{\ast}(G)$ induced by the conjugation action of $G$.

\medskip
It is part of our strategy (in the case of linear groups) to show, after dealing with the affine and symplectic cases of Section 4,  that we may restrict our calculations to groups with property ${\bf Q}$ which are primitive subgroups $G$ of ${\rm GL}(n,\mathbb{C})$ with $F(G) = Z(G)$, and $G/Z(G)$ almost simple, and where  $E = E(G)$ also acts (irreducibly and) primitively. The fact that $E(G)$ can be assumed to be primitive is not essential for the proof of Theorem 1.2, but it does reduce the number of resdual cases which need to be considered, and may help to improve the bound $60$ for some properties ${\bf Q}.$

\medskip
We will focus most attention on the set $A({\bf Q})$. The next result will be made more precise later, but it indicates that in some ways, the set $B({\bf Q})$ is redundant in the present context. The invariant $\beta({\bf Q})$ (when it exists) is, however, sometimes useful  to us for inductive purposes.

\medskip
\noindent {\bf Lemma 2.1:} \emph{ The set $A({\bf Q})$ is bounded above if and only if the set $C({\bf Q})$ is bounded above.
If either of these sets is bounded above, then the set $B({\bf Q})$ is bounded above.}

\medskip
\noindent {\bf Proof:} It is clear that $A({\bf Q})$ is bounded above if $C({\bf Q})$ is bounded above, since whenever $G$ is  a finite subgroup of ${\rm GL}(n,\mathbb{C})$ with property ${\bf Q},$  every Abelian normal subgroup of $G$ is contained in $F(G).$

\medskip
On the other hand, suppose that $A({\bf Q})$ is bounded above. We claim that $B({\bf Q})$ and $C({\bf Q})$ are both bounded above. To obtain an upper bound for $B({\bf Q}),$ it suffices to consider primitive permutation groups. If $G$ is isomorphic to a primitive subgroup of the symmetric group ${\rm Sym}_{n}$ with property ${\bf Q}$, then $G$ is certainly isomorphic to a subgroup of ${\rm GL}(n,\mathbb{C})$, so that $[G:F(G)] \leq \alpha({\bf Q})^{n-1},$ where $\alpha({\bf Q})$ is the least upper bound of $A({\bf Q}).$

\medskip
There are two cases to consider: either $F(G) = 1$, in which case we have  $|G| \leq \alpha({\bf Q})^{n-1},$ or else
$F(G) \neq 1$, in which case $G = F(G)H$, where $F(G)$ is an elementary Abelian normal $q$-group for some prime $q$ and $H$ is a maximal subgroup of $G$ with $F(G) \cap H = 1.$ In the second case, we have $|H| = [G:F(G)]$ and $n = q^{a}$ for some prime $q$ and some positive integer $a$. Then  $$|G| =n|H| \leq 2^{n-1}|H| \leq (2\alpha({\bf Q}))^{n-1}.$$ In either case, we have 
$$|G| \leq (2\alpha({\bf Q}))^{n-1}.$$ Hence $B({\bf Q})$ is bounded above by $2\alpha({\bf Q})$, so that the least upper bound 
$\beta({\bf Q})$ exists.

\medskip
Now, since $\beta({\bf Q})$ exists, the determination of whether $C({\bf Q})$ is bounded above may be reduced to considering primitive subgroups $G \leq {\rm GL}(n,\mathbb{C})$  with property ${\bf Q}.$ 

\medskip
Let $G$ be a finite primitive subgroup of  ${\rm GL}(n,\mathbb{C}).$ Then all Abelian normal subgroups of $G$ are central. 
If $F(G) = Z(G)$, then $$[G:Z(G)] = [G:F(G)] \leq \alpha({\bf Q})^{n-1}.$$ If $F(G) > Z(G)$, then all characteristic Abelian subgroups of $F(G)$ are central in $G$, so that the Fitting subgroup $F(G)$ is nilpotent of class $2.$ Then every element of $F(G) \backslash Z(G)$ has trace zero. A character calculation yields that\\ $[F(G):Z(G)] \leq n^{2}.$ Hence, in either case, we have $$[G:Z(G)] \leq n^{2}\alpha({\bf Q})^{n-1} \leq (4\alpha({\bf Q}))^{n-1}$$ so we deduce that $C({\bf Q})$ is bounded above by $4\alpha({\bf  Q}),$ and that its least upper bound 
$\gamma({\bf Q})$ exists, and does not exceed $4\alpha({\bf Q}).$

\medskip
\noindent {\bf Remark 2.2:} We will presently see that when ${\bf Q}$ is the property of solvability, then we have $\alpha({\bf Q}) = 6$ and 
$\gamma({\bf Q}) = 24,$ so that the equality $\gamma({\bf Q}) = 4\alpha({\bf Q})$ can sometimes be attained.

\medskip
We will show relatively easily that $\beta({\bf Q})$ does not exist if infinitely many alternating groups have property ${\bf Q},$ so that neither $\alpha({\bf Q})$ nor $\gamma({\bf Q})$ exists in that case. At this point, we will restrict our attention to 
properties ${\bf Q}$  which only a finite number of (non-isomorphic) alternating groups possess.

\medskip
We outline our proof of 1.2 now. We will show  that if $G$ is a primitive subgroup of a symmetric group
${\rm Sym}_{n}$ with $F(G) \neq 1$, then $$|G| \leq 24^{\frac{n-1}{3}}< 3^{n-1}.$$  The proof is independent of the choice of ${\bf Q}$, but it is possible to  improve this bound for some choices of ${\bf Q}.$  Similarly, we prove that if $G$ is a finite group with a non-central extraspecial normal $q$-subgroup $U$ of order $q^{2m+1}$ for some prime $q$, then $$[G:UC_{G}(U)]  \leq 720^{   \frac{q^{m}-1}{3}   } < 9^{q^{m}-1}$$ and $$[G:C_{G}(U)] \leq 24^{q^{m}-1 }.$$ This result is helpful as part of an inductive argument to obtain the correct type of bound for the order of primitive subgroups $G$ of ${\rm GL}(n,\mathbb{C})$ with property ${\bf Q}.$

\medskip
Again, this intermediate bound does not depend on ${\bf Q}$, but it may be possible to improve it for certain choices of ${\bf Q}.$

\medskip
Next, we use the classification of finite simple groups to prove that 

\medskip 
$\{ [G:Z(G)]^{\frac{1}{n-1}}: G \leq {\rm GL}(n,\mathbb{C})$ \emph{is primitive}, $ G$ \emph{has property} ${\bf Q}, $ \\$Z(G) = \Phi(G), G/Z(G)$ \emph{is almost simple}, $E(G)$ \emph{is primitive} $\}$  is bounded above. 

\medskip
\noindent {\bf Definition 2.3: } We let $\ell({\bf Q})$ denote the least upper bound of the set above (if the set is empty, we set $\ell({\bf Q}) = 1$).

\medskip
A substantial part of our proof of Theorem 1.2 is to prove by an inductive argument that $A({\bf Q})$ is bounded above by 
$${\rm max}\{720^{\frac{1}{3}}, \ell({\bf Q}) \}.$$  We also outline the modifications necessary to show that we likewise have 
$$\gamma({\bf Q}) \leq  {\rm max}\{24, \ell({\bf Q}) \},$$ and also $$\beta({\bf Q}) \leq {\rm max}\{720^{\frac{1}{3}}, \ell({\bf Q}) \}.$$ 

\medskip 
Hence whenever $G$ is a subgroup of ${\rm GL}(n,\mathbb{C})$ with property ${\bf Q}$, we have $$[G:F(G)] \leq  \left( {\rm max}\{\ell({\bf Q}), 720^{\frac{1}{3}} \} \right)^{n-1},$$  
$G$ has an Abelian normal subgroup of index at most $$ \left( {\rm max}\{ \ell({\bf Q}), 24 \} \right) ^{n-1},$$ 

\medskip
\noindent and whenever $G$ is a subgroup
of ${\rm Sym}_{n}$ with property ${\bf Q}$, we have $$|G| \leq  \left( {\rm max}\{ 720^{\frac{1}{3}}, \ell({\bf Q})  \}  \right)^{n-1}.$$  

\medskip
To deal with $\beta({\bf Q})$, we reduce to considering primitive permutation groups as usual. Now if $G$ 
is a primitive subgroup of ${\rm Sym}_{n}$ with property ${\bf Q}$, we have $$|G| \leq 24^{\frac{n-1}{3}}$$  if $F(G) \neq 1$. If $F(G)= 1$, then since $G$ is isomorphic to a subgroup of\\  ${\rm GL}(n-1,\mathbb{C}),$ we have $$|G| = [G:F(G)] \leq \alpha({\bf Q})^{n-2}.$$  In either case, we have $$|G| \leq \left({\rm max} \{\ell({\bf Q}), 720^{\frac{1}{3}} \} \right)^{n-1}.$$

\medskip
We return to our consideration of $\alpha({\bf Q}).$ It will be clear at this point that $\alpha({\bf Q})$ exists and that 
$$\ell({\bf Q}) \leq \alpha({\bf Q}) \leq {\rm max}\{ 720^{\frac{1}{3}}, \ell({\bf Q})  \}.$$ In particular, if $\ell({\bf Q}) \geq 720^{\frac{1}{3}},$
then $\alpha({\bf Q}) = \ell({\bf Q}).$ If $\ell({\bf Q}) < 720^{\frac{1}{3}},$ then $\alpha({\bf Q})$ exists and is less than or
equal to ${\rm 720}^{\frac{1}{3}} < 9.$  Hence, unless $\alpha({\bf Q})$ is already rather small, $\alpha({\bf Q})$ may be calculated by considering 
finite primitive subgroups $G$ of ${\rm GL}(n,\mathbb{C})$ with property ${\bf Q}$ such that $G/Z(G)$ is almost simple and $E(G)$ also acts (irreducibly and) primitively.

\medskip
To provide an explicit upper bound for $\ell({\bf Q}),$ we use the classification of the finite simple groups to prove that we have  $\ell({\bf Q})\leq 60$  if no non-Abelian simple alternating group ${\rm Alt}_{m}$ with $m >  151$ has property ${\bf Q}.$ 

\medskip
Since there are only finitely many sporadic simple groups, and since we are assuming that only finitely many alternating groups have property ${\bf Q}$, the existence of an upper bound of some kind reduces to consideration of finite simple groups of Lie type.

\medskip
We are able to prove that $$\left( \frac{m!}{2}\right)^{\frac{1}{m-2}}\leq \ell({\bf Q}) \leq {\rm max} \{60,m!^{\frac{1}{m-2}}\}$$ where $m$ is the maximal integer such that  ${\rm Alt}_{m}$ has property ${\bf Q}$.

\medskip
When $G$ is a finite primitive subgroup of ${\rm GL}(n,\mathbb{C})$ with property ${\bf Q}$ with $G/Z(G)$ almost simple, with $E = E(G)$ acting irreducibly, and with\\ $X = E/Z(E)$ a sporadic simple group, then we are able to prove that\\ $[G:Z(G)] < 15^{n-1}.$

\medskip
It requires a fairly detailed analysis to obtain $[G:Z(G)] \leq 60^{n-1}$ 
when $G$ is a finite primitive subgroup of ${\rm GL}(n,\mathbb{C})$ with property ${\bf Q}$ and with $G/Z(G)$ quasi-simple, with $E = E(G)$ acting irreducibly, such that $X = E/Z(E)$ is a finite simple group of Lie type.

\medskip
From now on, ${\rm s}({\bf Q})$ denotes the set of non-Abelian simple groups with property ${\bf Q},$ while ${\rm as}({\bf Q})$ denotes the set of almost simple groups with property ${\bf Q}$. Also, ${\rm pas}({\bf Q})$ denotes the set of groups $G$ with property ${\bf Q}$ such that $G/Z(G) \in {\rm as}({\bf Q}).$ 

\medskip
\section{Bounds imposed by the presence of alternating groups}

\medskip
\noindent {\bf Remark 3.1 :} It is well known and easy to check that the ratio $\frac{n!}{n^{n}}$ strictly decreases as $n$ increases. Hence the ratio  $\frac{n!^{\frac{1}{n}}}{n}$ certainly strictly decreases as $n$ increases (and always takes values in $(0,1]$). Thus the sequence $\left( \frac{n!^{\frac{1}{n}}}{n} \right)$ converges to its greatest lower bound. Since the limit of that sequence is well known to be $\frac{1}{e}$, we have $$n!^{\frac{1}{n}} > \frac{n}{e}$$ for every positive integer $n$, and certainly $$n!^{\frac{1}{n-1}} > \frac{n}{e}$$ for every integer $n \geq 5.$  Hence we have $$\left( \frac{n!}{2} \right)^{\frac{1}{n-1}} > \frac{n}{2e}$$ for every integer $n \geq 5.$ Using Lemma 2.1, we have: 

\medskip
\noindent {\bf Proposition 3.2:} \emph{Let ${\bf Q}$ be an extant, normally hereditary, property of finite groups. Suppose that ${\rm Alt}_{n} \in s({\bf Q})$ for some $n \geq 5$. Then $B({\bf Q})$ contains a real number greater than $\frac{n}{2e}.$ In particular, if $s({\bf Q})$ contains infinitely many alternating groups, then $B({\bf Q})$ is unbounded, and none of $\alpha({\bf Q}),\beta({\bf Q})$  or $\gamma({\bf Q})$ exists.}

\medskip
\noindent {\bf Remark 3.3:} For every $n >2$, we have $$\frac{n!}{n^{n-2}} > \frac{(n+1)!}{(n+1)^{n-1}},$$ so certainly $$\frac{ n!^{ \frac{1}{n-2} }  }{n} > \frac{(n+1)!^{\frac{1}{n-1}}}{n+1}.$$ 
Hence the sequence $$\left(\frac{n!^{\frac{1}{n-2}}}{n}\right)_{n \geq 3}$$ strictly decreases as $n$ increases, and also converges to $\frac{1}{e}.$

\medskip
In the other direction, for $n \geq 5$, 
we have $n! \leq  24n^{n-4} <  n^{n-2}$, so we always have $$ \frac{n}{e} < n! ^{\frac{1}{n-2}} < n.$$

\medskip
Now we also have $$\frac{ (n+1)!^{\frac{1}{n-1}}}{n!^{ \frac{1}{n-2}} }
=   \frac{(n+1)^{\frac{1}{n-1}} }{n!^{\frac{1}{(n-1)(n-2)}}} > \frac{(n+1)^{\frac{1}{n-1}} }{n^{\frac{1}{n-1}}} > 1,$$ so that the sequence 
$$\left(n!^{\frac{1}{n-2}}\right)_{n \geq 5}$$ increases strictly  as $n$ increases, as does the sequence $$\left(\left(\frac{n!}{2}\right)^{\frac{1}{n-2}}\right)_{n \geq 5}$$
(the last remark may be relevant when ${\rm Alt}_{n} \in {\rm s}({\bf Q})$, but ${\rm Sym}_{n} \not \in {\rm as}({\bf Q})).$

\medskip
\noindent {\bf Remark 3.4:} For $n > 7,$ the smallest degree of a non-trivial complex irreducible character of the Schur cover of ${\rm Alt}_{n}$ is $n-1$, as was known to I. Schur.\\ For $n =5,6$ the groups ${\rm Alt}_{5}$ and ${\rm Alt}_{6}$ (which may both be considered to be simple groups of Lie type) have respectively a double cover with  a complex representation of degree $2$ and a triple cover with a complex representation of degree $3$. Neither of these representations extends to a covering group of any larger subgroup of the automorphism group, as may be seen by considering the respective Sylow $5$-normalizers. Also, a double cover of ${\rm Alt}_{7}$ has a $4$-dimensional irreducible complex representation which does not extend to the covering group of  ${\rm Sym}_{7},$ on consideration of the Sylow $7$-normalizer. There is another exceptional case, corresponding to the $4$-dimensional representation of the double cover of ${\rm Alt}_{6}$ (which is isomorphic to ${\rm SL}(2,9)).$

\medskip		
\noindent {\bf Remark 3.5:} Careful use of Stirling's approximation shows that $151! < 60^{149},$ while $152! > 60^{150}.$ Hence we have 
$$n!^{\frac{1}{n-2}} < 60$$ if and only if $n < 152.$

\medskip
We now note:

\medskip
\noindent {\bf Lemma 3.6:} \emph{ If ${\rm SL}(2,5) $ has property ${\bf Q}$, then we have $\alpha({\bf Q}) \geq 60.$ If $3.{\rm Alt}_{6}$ has property ${\bf Q},$ then we have
 ${\rm \alpha}({\bf Q}) \geq \sqrt{360},$ and if $2.{\rm Alt}_{7} $ has property ${\bf Q}$, then we have 
$\alpha({\bf Q}) \geq 2520^{\frac{1}{3}}.$ If $m > 7$ is the largest integer $m$ such that ${\rm Alt}_{m} \in{\rm s}({\bf Q})$, 
then we have $$\alpha({\bf Q}) \geq \left( \frac{m!}{2} \right) ^{\frac{1}{m-2}}.$$  If ${\rm Sym}_{m} \in {\rm as}({\bf Q}),$ 
then we have $$\alpha({\bf Q}) \geq m!^{\frac{1}{m-2}}.$$ }

\medskip
\noindent {\bf Remark 3.7:} The lower bounds imposed on $\alpha({\bf Q})$ in Lemma 3.3 do not exceed $60$ unless $m > 151.$

\medskip
\section{ Affine, Symplectic  and extended Symplectic Cases}

\medskip
In the theory of finite permutation groups, the key residual objects to deal with after standard reductions are usually primitive permutation groups. If $G$ is a primitive subgroup of the symmetric group ${\rm Sym}_{n}$ (a transitive subgroup $G$ whose point-stabilizer $H$ is a maximal corefree subgroup of $G$), then if $F(G) \neq 1$, we must have 
$$n = |F(G)| = q^{m}$$ for some prime $q$ and positive integer $m$, while $H \cap F(G) = 1$ and $H$ is isomorphic to an irreducible subgroup of the general linear group 
${\rm GL}(m,q).$ On the other hand, if $F(G) = 1$, then $G = MH$ where $M$ is a minimal normal subgroup which is a direct product of isomorphic non-Abelian simple groups. We refer to the case $F(G) \neq 1$ as the \emph{affine case}.

\medskip
Similarly, in the theory of finite complex linear groups, the key objects to\\ study are primitive linear groups. If $G$ is a finite primitive subgroup of ${\rm GL}(n,\mathbb{C})$, and we consider a normal subgroup $M$ of $G$ minimal subject to strictly containing $Z(G)$ , there are two possibilities to consider. In the case that $F(G) > Z(G)$, it is possible to choose $M$ to be of the form $Z(G)Q$, where $Q$ is extra-special of order $q^{2m+1}$ for some prime $q$ and integer $m.$  If $q = 2$, then $Q$ necessarily has exponent $4$, and in the case that  $m = 1$ and $|Z(G)|$ is not divisible by $4$, $Q$ is quaternion of order $8$. We have  $Q = \Omega_{1}(O_{q}(M))$ if $q$ is odd and 
$Q = [\Omega_{2}(O_{2}(M)),G]$ if $q = 2$. Also, $q^{m}$ is a divisor of $n$, and $G/MC_{G}(M)$ is isomorphic to an irreducible subgroup of 
${\rm Sp}(2m,q)$ (by minimality of $M$). In our context, we are interested in bounding both $[G:Z(G)]$ and $[G:F(G)].$ It is clear that 
$[G:F(G)] \leq [G:MC_{G}(M)][C_{G}(M):F(C_{G}(M))]$ and that 
$[G:Z(G)] = [G:C_{G}(M)][C_{G}(M):Z(C_{G}(M))]$  (note that we have $Z(C_{G}(M)) = Z(G)$ in this situation, since all Abelian normal subgroups of $G$ are central by primitivity).

\medskip
To deal with the respective primitive cases, therefore, we will have an interest in bounding the quantities 
$${\rm ea}(q,m) = (q^{m}|{\rm GL}(m,q)|)^{\frac{1}{q^{m}-1}},$$ which we refer to as the \emph{extended affine bound}, $${\rm a}(q,m) = |{\rm GL}(m,q)|^{\frac{1}{q^{m}-1}},$$ which we refer to as the \emph{affine bound}, $${\rm es}(q,m) = (q^{2m}|{\rm  Sp}(2m,q)|)^{\frac{1}{q^{m}-1}},$$ which we refer to as the  \emph{extended symplectic bound} and $${\rm s}(q,m) = |{\rm Sp}(2m,q)|^{\frac{1}{q^{m}-1}},$$ which we refer to as the \emph{symplectic bound}. Bounding these quantities will ultimately allow us to turn our attention in the linear group case to the case that $F(G) = Z(G).$

\medskip
Notice that $${\rm a}(q,m) < q^{ \frac{m^{2}}{(q^{m}-1)}},$$  $${\rm ea}(q,m)  < q^{\frac{m^{2}+m} {(q^{m}-1)}},$$ 
$${\rm s}(q,m) < q^{    \frac{2m^{2}+m}{q^{m}-1}}$$  and $${\rm es}(q,m) < q^{\frac{2m^{2}+3m} {q^{m}-1}}.$$

\medskip
We consider the behaviour of the functions $$f(x,y) = f_{(a,b)}(x,y) = \exp{ \left( \frac{(\log{y}) (ax^{2}+bx)}{y^{x}-1} \right)},$$ where $a,b$ are specified non-negative real numbers, $x \geq 1, y >  1.$ Notice that for every $x$, the function $f_{x}: f_{x}(y) = f(x,y)$ increases or decreases as $y$ increases in the same manner as 
$$g_{x}: g_{x}(y) =  \frac{\log{y}}{y^{x}-1}$$ does.

\medskip
Now $$g_{x}^{\prime}(y) = \frac{1}{y(y^{x}-1)} - \frac{(x\log{y})y^{x-1}}{(y^{x}-1)^{2}}$$ is negative whenever $$ \log{y} > \frac{y^{x}-1}{xy^{x}},$$ so certainly whenever $y > e^{\frac{1}{x}}.$ In other words, $g_{x}^{\prime}(y)< 0 $ whenever $y^{x} > e$.

\medskip
Thus for any fixed $x \geq 1$, the maximum value of $f_{(a,b)}(x,y)$ occurs for $$y \in (1,e^{\frac{1}{x}})$$ and $f_{(a,b)}(x,y)$ decreases monotonically for $y \in (t,\infty)$, where $t$ is the least integer greater than $e^{\frac{1}{x}}.$

\medskip
This means that the prime $q$ for which $f_{(a,b)}(1,q)$ is maximal is either $q = 2$ or $q = 3$, while for any integer $d > 1$, the prime $q$ for which $f_{(a,b)}(d,q)$ is maximal is $q = 2.$

\medskip
Now we can prove : 

\medskip
\noindent {\bf Proposition 4.1:}  \emph{ For all primes $q$ and positive integers $m$, we have: $${\rm es}(q,m)  \leq {\rm es}(2,1) = 24,$$
$${\rm s}(q,m)  \leq {\rm s}(2,2) = (720)^{\frac{1}{3}} < 9,$$ 
$${\rm ea}(q,m)  \leq {\rm ea}(2,2) = 24^{\frac{1}{3}},$$ 
$${\rm a}(q,m)  \leq {\rm a}(2,3) = 168^{\frac{1}{7}}.$$ }

\medskip
\noindent {\bf Proof:} We will make use of the fact that (by the binomial theorem), we have $2^{m}-1 \geq \frac{m^{2}+m}{2}$ for  each integer $m \geq 2$ and $2^{m}-1 \geq \frac{m^{3}+5m}{6}$ for each integer $m \geq 3.$ Also, we have $3^{m}- 1 \geq 2m^{2}$ for each integer $m \geq 2.$

\medskip
We see that $${\rm es}(2,1) = 4|{\rm SL}(2,2)| = 24,$$ $${\rm es}(3,1) \leq f_{2,3}(1,3) = 3^{\frac{5}{2}} < 16.$$
$${\rm es}(2,2) = (16|{\rm Sp}(4,2)|)^{\frac{1}{3}}< 4 \times (180)^{\frac{1}{3}} < 24.$$
$${\rm es}(3,2) \leq f_{2,3}(2,3) = 3^{\frac{14}{8}} < 9.$$

\medskip
For $q >3$, we have $${\rm es}(q,1) \leq f_{2,3}(1,q) \leq f_{2,3}(1,3) < 16$$ and also $${\rm es}(q,2) \leq f_{2,3}(2,q) \leq f_{2,3}(2,3) < 9.$$

\medskip
For $m \geq 3$, we have $${\rm es}(q,m) \leq f_{2,3}(m,q) \leq f_{2,3}(m,2) \leq 2^{\frac{6m^{2}+9m}{2m^{2}+m}} < 16.$$ Thus we have ${\rm es}(q,m) \leq 24$
in all cases.

\medskip
We see that $${\rm s}(2,1) = |{\rm SL}(2,2)| = 6,$$ $${\rm s}(3,1) \leq f_{2,1}(1,3) = 3^{\frac{3}{2}} < 6.$$
$${\rm s}(2,2) = (|{\rm Sp}(4,2)|)^{\frac{1}{3}}< 9.$$
$${\rm s}(3,2) \leq f_{2,1}(2,3) = 3^{\frac{5}{4}} < 4.$$

\medskip
For $q >3$, we have $${\rm s}(q,1) \leq f_{2,1}(1,q) \leq f_{2,1}(1,3) < 6$$ and also $${\rm s}(q,2) \leq f_{2,1}(2,q) \leq f_{2,1}(2,3) < 4.$$

\medskip
For $m \geq 3$, we have $${\rm s}(q,m) \leq f_{2,1}(m,q) \leq f_{2,1}(m,2) \leq 2^{\frac{6m^{2}+3m}{2m^{2}+m}} = 8.$$ Thus we have 
$${\rm s}(q,m) \leq (720)^{\frac{1}{3}}$$ in all cases.

\medskip
We see that $${\rm ea}(2,1) = 2,$$ $${\rm ea}(3,1) = \sqrt{6},$$
$${\rm ea}(2,2) = (24)^{\frac{1}{3}},$$    $${\rm ea}(2,3) = (8 \times 168)^{\frac{1}{7}} < {\rm ea}(2,2),$$
$${\rm ea}(3,2) = (9 \times 48)^{\frac{1}{8}} < 21^{\frac{1}{4}} < {\rm ea}(2,2).$$

\medskip
For $q  \geq 5$, we have $${\rm ea}(q,1) = (q(q-1))^{\frac{1}{q-1}} \leq 20^{\frac{1}{4}} < {\rm ea}(2,2). $$ 

\medskip
For $m \geq 3, q \geq 3$, we have $${\rm ea}(q,m) \leq f_{1,1}(m,q) \leq f_{1,1}(m,3) \leq 3^{\frac{m^{2}+m}{2m^{2}}}  \leq 3^{\frac{2}{3}} < {\rm ea}(2,2).$$ 

\medskip
For integers $m \geq 4$, we see by induction that $$\frac{m^{2}+m}{2^{m}-1} \leq \frac{4}{3},$$ and so for all primes $q$, we have 
$${\rm ea}(q,m) \leq f_{1,1}(m,q) \leq f_{1,1}(m,2)  \leq 16^{\frac{1}{3}} < {\rm ea}(2,2).$$
	
\medskip
Thus we have $${\rm ea}(q,m) \leq (24)^{\frac{1}{3}}$$
in all cases.

\medskip
We see that $${\rm a}(2,1) = 1,$$ $${\rm a}(3,1) = {\rm a}(5,1) = \sqrt{2},$$
$${\rm a}(2,2) = 6^{\frac{1}{3}} < 2,$$ $${\rm a}(2,3) = 168^{\frac{1}{7}},$$ $${\rm a}(3,2) = 48^{\frac{1}{7}} < {\rm a}(2,3).$$.

\medskip
For $q  \geq 3$, we have $${\rm a}(q,1) \leq f_{1,0}(1,q) \leq f_{1,0}(1,3) = \sqrt{3} < 2.$$ 

\medskip
For $q  >  3$ and $m \geq 2$ we have $${\rm a}(q,m) \leq f_{1,0}(m,q) \leq f_{1,0}(m,3) = 3^{\frac{m^{2}}{3^{m}-1}} < \sqrt{3}.$$ 

\medskip
When $m \geq  3$ is an integer, we may check by induction that  $\frac{m^{2}}{2^{m}-1}$ decreases as $m$ increases. When $m = 4$, we see that $$|{\rm GL}(4,2)| < 2^{15},$$
so that ${\rm a}(2,4) < 4.$ When $m \geq 5$ and $q$ is a prime, we have $${\rm a}(q,m) \leq f_{1,0}(m,q) \leq f_{1,0}(m,2)  \leq 2^{\frac{25}{31}} < 2.$$

\medskip
Hence in all cases, we have $${\rm a}(q,m) \leq {\rm a}(2,3) = 168^{\frac{1}{7}}.$$

\medskip
The relevance of these results to the questions we are considering is exhibited now:

\medskip
\noindent {\bf Corollary 4.2:} \emph{ a)  Let $G$ be a primitive permutation group of prime power degree $d = q^{m}$ with $F(G) \neq 1$. Then 
$F(G)$ is the unique largest Abelian normal subgroup of $G$. Furthermore, $$|G| \leq  24^{\frac{d-1}{3}} < 3^{d-1}$$ and 
$$[G:F(G)] \leq  168^{\frac{d-1}{7}}.$$ }

\medskip
\noindent  \emph{b)  Let $G$ be a finite primitive subgroup of ${\rm GL}(n, \mathbb{C})$. Suppose that $F(G) \neq Z(G)$, and let $M$ be a normal subgroup of $G$ with $Z(G) < M \leq F(G)$, and minimal subject to that. Then there is a prime power $q^{m}$ which divides $n$ such that $M = Z(G)Q,$ where $Q$ is extraspecial of order $q^{2m+1}$  and, furthermore, $$[G:MC_{G}(M)] \leq 720^{\frac{q^{m}-1}{3}} < 9^{q^{m}-1}.$$  Also, $$[G:C_{G}(M)] \leq 24^{q^{m}-1}.$$ }

\medskip
\noindent {\bf Remark 4.3:} If our property ${\bf Q}$ is extant and irredundant, then the inequalities of 
Corollary 4.1  still hold, but may sometimes be improved. For example, when all groups which have property ${\bf Q}$ are solvable, the bound $720^{\frac{1}{3}}$ can be improved to $6$ (or sometimes less, for example the bound can be improved to $\sqrt{6}$ if, in addition, all Sylow subgroups are Abelian). Similarly, the bound $168^{\frac{d-1}{7}}$ can be improved if ${\rm GL}(3,2) \not \in {\rm s}({\bf Q}).$

\medskip
\section{Some invariants of non-Abelian simple groups}

\medskip
Let $X$ be a non-Abelian finite simple group, and let $L$ be the (unique up to isomorphism) maximal perfect central extension of $X$, which we have already referred to as the Schur cover of $X$. Let $\chi$ be a non-trivial  (not necessarily faithful) irreducible character of $L$. Let $I(\chi)$ denote the inertial subgroup of $\chi$ in $C_{{\rm Aut}(L)}(Z(L))$ (this last subgroup is isomorphic to  a subgroup of ${\rm Aut}(X)$ and $I(\chi)/L$ is isomorphic to a subgroup of ${\rm Out}(X)).$ 

\medskip
\noindent {\bf Definition 5.1:} We define the \emph{maximum  projective ratio} of $X$, denoted ${\rm mpr}(X)$,  to be $$ {\rm max}_{\{1 \neq \chi \in {\rm Irr(L)}\}} (|X|[I(\chi):L])^{\frac{1}{(\chi(1)-1)}}.$$
While the inertial subgroups $I(\chi)$ may be difficult to determine exactly, it is clear that $$ |X|^{\frac{1}{d-1}} \leq {\rm mpr}(X) \leq |{\rm Aut}(X)|^{\frac{1}{d-1}},$$ where $d$ is the minimal degree of a  non-trivial complex  irreducible  character of $L$. 

\medskip
\noindent {\bf Remark 5.2:} The motivation for the definition of ${\rm mpr}(X)$ is that whenever $G$ is a finite primitive subgroup of 
${\rm GL}(n,\mathbb{C})$ such that $G/Z(G)$ is almost simple with $E = E(G)$ irreducible and $E/Z(E) \cong X$, then we have 
$$[G:Z(G)]^{\frac{1}{n-1}} \leq {\rm mpr}(X),$$ since $G/EZ(G)$ is isomorphic to a subgroup of ${\rm Out}_{G}(X)$ which stabilizes the irreducible character of $E$ afforded by the given representation.

\medskip
Thus, returning to our property ${\bf Q},$ the least upper bound $\ell({\bf Q})$ exists and is bounded above by the least upper bound 
of $\{ {\rm mpr}(X) : X \in {\rm s}({\bf Q}) \}$ if this last set is bounded above (in the case that ${\rm s}({\bf Q})$ is non-empty). 

\medskip
We also recall that every finite non-Abelian simple group $X$ may be generated by two elements, (a fact which presently requires the Classification of the Finite Simple Groups), so that we always have  $|{\rm Aut}(X)| < |X|^{2}.$ Hence we have:

\medskip
\noindent {\bf Lemma 5.3:} \emph{ For any non-Abelian finite simple group $X$, we have 
 $${\rm mpr}(X) < |X|^{\frac{2}{d-1}},$$ where $d$ is the minimum degree of a non-trivial complex irreducible character of the Schur cover of $X.$}

\medskip
\section{ Bounding the maximum projective ratio for the sporadic finite simple groups}

\medskip
\noindent {\bf Remark 6.1:} The minimal degrees of non-trivial irreducible characters of Schur covers of  the sporadic simple groups are conveniently listed in the paper [7] of C. Jansen, and can also be found in [5]. Since each sporadic simple group has outer automorphism group of order at most two, it also follows that when $X$ is a sporadic finite simple group, we have $${\rm mpr}(X) \leq (2|X|)^{\frac{1}{d-1}},$$ where $d$ is the minimal degree of a non-trivial irreducible complex character of the Schur cover of $X$. There are only two sporadic simple groups $X$ which have  ${\rm mpr}(X) >7$: these are the sporadic Suzuki group $Suz$, where we have $9 < {\rm mpr}(Suz) < 9.1$ (in this case, $d=12$) and Janko's second group $J_{2}$
(also known as the Hall-Janko group $HJ$), which has $|J_{2}| = 604800$ and $d=6$. We have $14 < {\rm mpr}(J_{2}) < 14.333$. In the last case, we use the fact that the irreducible character of degree $6$ of the Schur cover is not invariant under the full automorphism group, which can be seen by considering the Sylow $7$-normalizer.

\medskip
Hence we have:

\medskip
\noindent {\bf Proposition 6.2:} \emph{ Let $X$ be a finite sporadic simple group. Then we have ${\rm mpr}(X) < 15.$}

\medskip
\section{Bounding the maximum projective ratio for finite simple groups of Lie type and for alternating groups}

\medskip
The arguments of Section 4 may be adapted to produce explicit bounds on ${\rm mpr}(X)$ for simple groups of Lie type, using bounds such as those provided  by Landazuri and Seitz [8] and Tiep [11], on the minimal degree of non-trivial irreducible characters of the Schur covers of these groups. The bounds of Landazuri and Seitz have also previously been used in post-classification bounds obtained for Jordan's theorem on complex linear groups by M.J.Collins [3] and B. Weisfeiler [12,13], and Collins also made use of the bounds of Tiep [11]. Also, Collins tabulates some of the bounds from [8,11] in [3].

\medskip
The key point for most of the following calculations is that when $a$ is a non-negative real number, the function mapping $x$ to $x^{\frac{1}{x-a}}$ strictly decreases as $x$ increases for $x > {\rm max}(a,e)$. As in Section $3$, when trying to estimate the maximum value of this function in case $x$ is an integer greater than $1$,we usually only need to check the cases $x = 2$ and/or $x=3$ to locate the maximum value, although for some choices of $a$, we may need to consider larger integers.

\medskip
Our main aim in this section is to prove: 

\medskip
\noindent {\bf Proposition 7.1:} \emph{Let $X$ be  a finite simple group of Lie type. Then\\ ${\rm mpr}(X) \leq 60.$}

\medskip
\noindent {\bf Proof:}  We know that ${\rm mpr}(X) \leq |{\rm Aut}(X)|^{\frac{1}{d-1}}$, where $d$ is the minimum degree of a non-trivial irreducible character of the Schur cover of $X$.  The Landazuri-Seitz and/or Tiep bounds referenced above give a lower bound for $d$. In most cases below, $|{\rm Out}(X)|^{\frac{1}{d-1}}$ is small enough to be negligible in this context, and it suffices to obtain upper bounds for $|X|^{\frac{1}{d-1}}$ when $X$ is a finite simple group of Lie type, and $d$ is the minimum degree of a non-trivial  irreducible complex character of the Schur cover of $X$. Occasionally, we will need to replace $d$ by an estimate 
$d^{\prime} < d,$ which means that our estimate $|X|^{\frac{1}{d^{\prime} -1}}$  is a slightly larger upper bound than the true value of  $|X|^{\frac{1}{d-1}}.$ 

\medskip
We first consider $X = {\rm PSL}(n,q) $ where $n \geq 3$ . Here, we  may take \\ $d^{\prime} = q^{n-1} + 1$ as a lower bound for the minimum degree $d$, except for the case 
$q=3, n = 4$, where we take $d^{\prime} = 26$.  Then whenever $(q,n) \neq (3,4)$, we have 
$$|X|^{ \frac{1}{d^{\prime}-1}} <  \left( \left( q^{n-1} \right)^{ \frac{1}{q^{n-1}} }  \right)^{  \frac{ n^{2}}{n-1} } .$$ 

\medskip
For any fixed value of $n$, using the discussion above, we see that $q^{ \frac{n^{2}}{ q^{n-1}} }$ takes its maximum value (over integers $q >1$) when $q = 2$, when the value is  
$$2^{ \frac{n^{2}}{2^{n-1} }   },$$ which decreases as $n \geq 3$ increases, and  is always at most $2^{\frac{9}{4}}.$

\medskip
When $q = 3, n = 4,$ we obtain $$|X|^{\frac{1}{d^{\prime} -1}} < 3^{\frac{16}{25}}<  2.011.$$

\medskip
Let us examine $|X|^{\frac{1}{d-1}} $ when $X = {\rm PSL}(2,q)$ . Consider first the case that $q = 2^{m}$. We saw in the alternating groups discussion that 
${\rm mpr}({\rm SL}(2,4)) = 60.$ When $q = 2^{m} \geq 8$, we have $d = q-1$ and  $$|X|^{\frac{1}{d-1}} \leq q^{\frac{3}{q-2}} \leq \sqrt{8}.$$

\medskip
When $q$ is odd, we only need consider $q >5.$ In this case, $d = \frac{q-1}{2}$, and we have $$|X|^{\frac{1}{d-1}}< q^{\frac{6}{q-3}} < 7^{\frac{3}{2}} < 19.$$

\medskip
We may conclude from this discussion that ${\rm mpr}(X) \leq 60$ for each non-Abelian simple group $X \cong {\rm PSL}(n,q).$

\medskip
Now let us consider ${\rm PSp}(2n,q)$ when $n \geq 2$, the case $n = 1$ being covered by the ${\rm PSL}(2,q)$ case. We also note that we have dealt with $X = {\rm PSp}(4,2)^{\prime}$, which is isomorphic to ${\rm Alt}_{6}.$

\medskip
For $q$ odd, we have done most of the work already in section 4: the minimal non-trivial irreducible character degree of the Schur cover is 
$d = \frac{q^{n}-1}{2}$, and we obtain  $d-1 = \frac{q^{n}-3}{2}$ in that case. We obtain a bound $$q^{\frac{4n^{2}+2n}{q^{n}-3}}.$$ By the discussion at the beginning of the section (with $x = q^{n})$, we see that for any particular $n$, this expression is maximized when $q = 3$. Also, when $q=3$, we see that $3^{n}-3 \geq 2n^{2}-2$ for $n \geq 2$ by the binomial theorem.
Hence whenever $n \geq 2,$ we have $$3^{ \frac{4n^{2}+2n}{3^{n}-3} } \leq  3^{2 + \frac{n+2}{n^{2}-1} } \leq 3^{\frac{10}{3}}.$$ 

\medskip
When $q$ is even, a lower bound for the minimum degree $d$  is $$d^{\prime}= \frac{q^{n-1}(q^{n-1}-1)(q-1)}{2}.$$ If $q = 2,$ we may, and do, suppose that $n \geq 3$.
Then in all cases we have $d^{\prime} \geq (q^{n}-2)$ , and we obtain an upper bound $$|X|^{\frac{1}{d^{\prime}-1}} \leq q^{\frac{2n^{2}+n}{q^{n}-3}}.$$
For any given choice of $n> 2$, this is maximized  when $q = 2,$ where the value attained is at most $2^{\frac{21}{5}} < 20$. When $n =2$ and $q \geq 4$, we have the better bound $d^{\prime} = 2(q-1)^{2}$, so that $d^{\prime}-1  > 2q^{2}-4q$ and 
 $$ q^{\frac{10}{2q^{2}-4q}}$$ is maximal when  $q = 4,$ so we obtain $$|X|^{\frac{1}{d^{\prime}-1}} < 4$$ in all cases.

\medskip
For  $X = P\Omega(2n+1,q)$, $q$ odd, $n \geq 3$, we obtain the same group order as for ${\rm PSp}(2n,q),$ but in each case the lower bound $d^{\prime}$ is at least $\frac{q^{n}-1}{2}$, so that $|X|^{\frac{1}{d^{\prime}-1}}$ is at most as large as for ${\rm PSp}(2n,q).$ 

\medskip
For  $X = P\Omega^{+}(2n,q)$, $n \geq 4$, we obtain a slightly smaller group order than ${\rm PSp}(2n,q),$ but in each case (except $q=2, n = 4$) the lower bound $d^{\prime}-1$ is at least as large as the lower bound we used above for $X = {\rm PSp}(2n,q),$ so $$|X|^{\frac{1}{d^{\prime}-1}}$$ reduces in each case other than the exceptional one. When $n = 4, q = 2$, there is a complex irreducible representation of (the minimal possible) degree $8$ of the Schur cover, and direct computation shows that 
$|X|^{\frac{1}{7}} < 17.$

\medskip
If we compare  $X = P\Omega^{-}(2n,q)$, $n \geq 4$, with the symplectic case, we obtain a smaller group order than for ${\rm PSp}(2n,q).$ On the other hand, the lower bound $d^{\prime}$ is always larger than the lower bound we used above for ${\rm PSp}(2n,q)$
(when $q$ is odd, we have $d^{\prime } \geq   q(q^{n}+1) -1$, compared to $\frac{q^{n}-1}{2}$, and when $q$ is even, we have $$d^{\prime} \geq q^{n+1} + q -1$$
compared to the bound $d^{\prime} \geq q^{n}-2$ used above). Hence the upper bound $$|X|^{\frac{1}{d^{\prime}-1}}$$ is smaller than that obtained in the symplectic group case.

\medskip
Neglecting unitary groups for the moment, it is easy to check (using the Landazuri-Seitz and/or Tiep bounds), that if $X$ is one of the remaining exceptional or twisted type Lie groups, we have $|X|^{\frac{1}{d-1}} < 10.$ In particular, for the Suzuki group $Sz(2^{2m+1}),$ we have $d = 2^{m}(2^{2m+1}-1)$ and  $$|Sz(2^{2m+1})|^{\frac{1}{d-1}} < 4\sqrt{2}.$$

\medskip
To deal with unitary groups, we may note that (pairing successive terms), we have $$\prod_{i=1}^{n-1} (q^{i+1} - (-1)^{i+1}) < q^{\frac{n(n+1)}{2}},$$ so we have an upper bound of $q^{n^{2}+n}$ for $|{\rm U}(n,q)|$ (in fact, $q^{n^{2}+n-1}$, but it suffices to use the cruder estimate here).

\medskip
This time, the lower bound for the minimum degree $d$ is $$d^{\prime} = \frac{q^{n}- 1}{q+1} $$ if $n \geq 4$ is even, and is 
$$d^{\prime} = \frac{q^{n}- q}{q+1} $$ if $n \geq 3$ is odd.

\medskip
Then $\frac{1}{d^{\prime}-1}$  may be written as $$ \frac{ 1 + \frac{1}{q}}{q^{n-1}-1-\frac{2}{q}}$$ if $n$ is even, which is less than or equal to $$\frac{3}{2(q^{n-1} - 2)}.$$

\medskip
As before, for any particular choice of even $n \geq 4 $, $$q^{\frac{3n^{2}+3n}{2(q^{n-1} - 2)}}$$ decreases as $q \geq  2$ increases, so the maximum value is attained at $q = 2$, where it is always at most $2^{5}$, so we have $|X|^{\frac{1}{d^{\prime}-1}} \leq 32$ when $n \geq 4$ is even.

\medskip
For $n \geq 3$ odd,  the lower bound for the minimum degree $d$ is $$d^{\prime} = \frac{q^{n}- q}{q+1}, $$ but we note that ${\rm U}(3,2)$ is not simple and should not be considered.

\medskip
Then $\frac{1}{d^{\prime} -1}$  may be written as $$ \frac{ 1 + \frac{1}{q}}{q^{n-1}-2-\frac{1}{q}}$$ if $n$ is odd, which is less than or equal to $$\frac{3}{(2q^{n-1} - 5)}.$$

\medskip
As in earlier discussions,  $$q^{\frac{3n^{2}+3n}{(2q^{n-1} - 5)}}$$ decreases as  $q \geq 2$ increases and the maximum value is taken when $q = 2$ and $n = 5.$ We thus obtain  $$|X|^{\frac{1}{d^{\prime}-1}} \leq 2^{\frac{10}{3}} < 10.079.$$

\medskip
We need to consider the case $n =3, q \geq 3$ separately. In this case, $\frac{1}{d^{\prime}-1}$  may be written as $$ \frac{ 1 + \frac{1}{q}}{q^{n-1}-2-\frac{1}{q}}$$  which is less than or equal to $$\frac{4}{(3q^{n-1} - 7)}.$$ Then our upper bound for $|X|^{\frac{1}{d-1}}$ is maximized when $q  = 3$, where the value is  $3^{\frac{12}{5}}. $

\medskip
For completeness, we record the following consequence (in the present context) of the results of Section 3 (cf. Lemma 3.3):

\medskip
\noindent {\bf Proposition 7.2:} \emph{We have ${\rm mpr}({\rm Alt}_{m}) \leq 60$ if $m \leq 151$ and \\ 
${\rm mpr}({\rm Alt}_{m}) = m!^{\frac{1}{m-2}}$ if $m >  151.$}

\medskip
Recalling that ${\rm mpr}(X) \leq |X|^{\frac{2}{d-1}}$ for a finite simple group $X$, which already deals with many cases, Proposition 7.1 and 7.2
(and the classification of finite simple groups) together yield:

\medskip
\noindent {\bf Theorem 7.3:} \emph{ Let $X$ be a finite simple group with ${\rm mpr}(X) > 60.$ Then $X$ is an alternating group ${\rm Alt}_{m}$ with $m > 151.$}

\medskip
\emph{In particular, for our property ${\bf Q}$, we know that only finitely many simple alternating groups have property ${\bf Q},$
so $\{ {\rm mpr}(X) : X \in  {\rm s}({\bf Q}) \}$ is bounded above and hence the least upper bound $\ell({\bf Q})$ exists.}

\medskip
\noindent {\bf Remark 7.4:} More precisely,  it follows from Theorem 7.3 that  
$$\{ {\rm mpr}(X) : X \in  {\rm s}({\bf Q}) \}$$ is bounded above by $60$ when there is no alternating group ${\rm Alt}_{m}$ in ${\rm s}({\bf Q})$ with $m >151,$ and is bounded above by $m!^{\frac{1}{m-2}}$ when $m > 151$ is the largest integer with ${\rm Alt}_{m} \in {\rm s}({\bf Q}).$ 

\medskip
Hence the least upper bound $\ell({\bf Q})$ exists for our property ${\bf Q}$, and we have ${\ell}({\bf Q}) \leq {\rm max}\{ 60,(m!)^{\frac{1}{m-2}}\},$ where $m$ is the largest integer such that the non-Abelian simple group ${\rm Alt}_{m} \in {\rm s}({\bf Q}).$

\medskip
\section{Conclusion of the proof of Theorem 1.2}

\medskip
\noindent {\bf Summary  8.1:} We know now that the least upper bound $\ell({\bf Q})$ exists for our specified property ${\bf Q},$ and that $\ell({\bf Q}) \leq 60$ if ${\rm s}({\bf Q})$ contains no alternating group ${\rm Alt}_{m}$ with $m > 151$, while $\ell({\bf Q}) \leq (m!)^{\frac{1}{m-2}}$ if $m > 151$ is the largest integer such that ${\rm s}({\bf Q})$ contains the alternating group ${\rm Alt}_{m}.$ 

\medskip
We now prove a general Lemma which may be useful in other contexts.

\medskip
\noindent {\bf Lemma 8.2: } \emph{Let $H$ be a finite primitive subgroup of ${\rm GL}(d,\mathbb{C})$ for some integer $d > 1$, 
and let $U > Z(H)$ be normal subgroup of $H.$ Then $U$ is isomorphic to a subgroup of ${\rm GL}(e,\mathbb{C})$ where $e$ is a divisor of $d$,
and $C_{G}(U)$ is isomorphic to a subgroup of ${\rm GL}(h,\mathbb{C})$ where $h$ is a divisor of $f = \frac{d}{e}$. Furthermore,
$Z(UC_{H}(U)) = Z(H)$, and $H/UC_{H}(U)$ is isomorphic to a subgroup of ${\rm Out}_{H}(U),$ the group of outer automorphisms of $U$ induced by the conjugation action of $H$. Also, $Z(U) = Z(C_{H}(U)) = Z(H)$ and $F(UC_{H}(U)) = F(U)F(C_{H}(U)).$ Finally, we have both 
$$[H:Z(H)] \leq |{\rm Out}_{H}(U)|[U:Z(U)][C_{H}(U): Z(C_{H}(U))]$$ and $$[H:F(H)] \leq [H: F(UC_{H}(U))] \leq |{\rm Out}_{H}(U)|[U:F(U)][C_{H}(U): F(C_{H}(U))].$$}

\medskip
\noindent {\bf Proof:} The first statement follows by Clifford's theorem and primitivity, taking $e = \mu(1)$ for $\mu$ an irreducible constituent 
of ${\rm Res}_{U}^{G}(\chi)$, where ${\rm \chi}$ is the irreducible character defined by $\chi(g) = {\rm trace}(g)$ for all $g \in G.$ 
Note that $e > 1$ here.

\medskip
For the second statement, note that (by primitivity of $H$), the centralizer algebra $C_{M_{d}(\mathbb{C})}(U)$ is isomorphic to $M_{f}(\mathbb{C})$, where $d = ef,$  and that the group of units of this centralizer algebra is isomorphic to ${\rm GL}(f,\mathbb{C}).$ Clearly, $C_{H}(U)$ is isomorphic to a subgroup of this group of units, and we see that ${\rm Res}_{C_{H}(U)}^{H}(\chi)$ decomposes as $e$ times the character afforded by the embedding of $C_{H}(U)$ in ${\rm GL}(d,f).$ On the other hand, $C_{H}(U) \lhd H$ so that, by primitivity of $H$, ${\rm Res}_{C_{H}(U)}^{H}(\chi)$
decomposes as a multiple of a single irreducible character, of degree $h$, say. Then $h$ divides $f = \frac{d}{e}.$

\medskip
It is clear that $C_{H}(UC_{H}(U)) \leq UC_{H}(U),$ and hence that $$C_{H}(UC_{H}(U)) = Z(UC_{H}(U)) \lhd H.$$ Thus we have
$C_{H}(UC_{H}(U)) = Z(H)$ by primitivity. Now we have  $$Z(U) = U \cap C_{H}(U) = Z(C_{H}(U)) = Z(H).$$ Hence $UC_{H}(U)/Z(H)$ is isomorphic to  $U /Z(H) \times C_{H}(U)/Z(H).$ From this, we easily deduce that $F(UC_{H}(U) = F(U)F(C_{H}(U))$ and that $$F(H) \cap F(C_{H}(U)) = Z(H).$$
Thus $$[H:F(H)] \leq [H:UC_{H}(U)] [UC_{H}(U):F(UC_{H}(U))] \leq |{\rm Out}_{H}(U)|[U:F(U)][C_{H}(U):F(C_{H}(U))],$$ and we also have 
$$[H:Z(H)] \leq [H: Z(UC_{H}(U))] \leq |{\rm Out}_{H}(U)|[U:Z(U)][C_{H}(U): Z(C_{H}(U))].$$

\medskip
Now we complete the proof of Theorem 1.2 by virtue of the following Proposition, which is stronger than we need for the immediate purpose.

\medskip
\noindent{\bf Proposition 8.3:} \emph{ We have $\alpha({\bf Q}) \leq {\rm max}\{720^{\frac{1}{3}}, \ell({\bf Q}) \}.$ }

\medskip
\noindent {\bf Proof:} If the proposition is false, then we may choose a real number $$r > {\rm max}\{720^{\frac{1}{3}}, \ell({\bf Q})\}$$ 
such that there is a group $G$ with property  ${\bf Q}$ which is a finite subgroup of ${\rm GL}(n,\mathbb{C})$ with $[G:F(G)] \geq  r^{n-1}$. We may choose the pair $(G,n)$ with first $n$, then $[G:F(G)]$, then $[G:Z(G)]$, then $|G|$, minimized, subject to\\ $[G:F(G)] \geq r^{n-1}.$ Notice that we have $n > 1$ and $G \neq F(G).$ Also, we certainly have $r > 8.$ 

\medskip
Our goal is to prove that $G$ is primitive with $G/Z(G)$ is almost simple, $Z(G) = \Phi(G),$ and $E = E(G)$ acts primitively in the given representation.

\medskip
However, in that case, we have $$[G:Z(G)]^{\frac{1}{n-1}} \leq \ell({\bf Q})$$ by definition of the least upper bound $\ell({\bf Q}).$
Then certainly we have  $$[G:F(G)] = [G:Z(G)] \leq \ell({\bf Q})^{n-1} < r^{n-1},$$ contrary to assumption.

\medskip
Thus we may conclude that $$[G:F(G)] \leq \left( {\rm max} \{720^{\frac{1}{3}}, \ell({\bf Q})\} \right)^{n-1} $$ whenever $G$ is a subgroup of ${\rm GL}(n,\mathbb{C})$ with property ${\bf Q}.$ Hence we do, after all, have $$\alpha({\bf Q}) \leq  {\rm max}\{720^{\frac{1}{3}},\ell({\bf Q}) \}.$$

\medskip
By Remark 6.3, we may then conclude that $\alpha({\bf Q}) \leq 60$ when ${\bf Q}$ contains no alternating group ${\rm Alt}_{m}$ with 
$m > 151$, and that $\alpha({\bf Q}) \leq m!^{\frac{1}{m-2}}$ when $m > 151$ is the largest integer such that the alternating group ${\rm Alt}_{m}$ has property ${\bf Q}$. 

\medskip
Let us return to the proof proper. By the minimality of the pair $(G,n)$, we see that $G$ is primitive as linear group. This requires us to note that if $H$ is isomorphic to a subgroup of the symmetric group ${\rm Sym}_{t}$ for some $t \leq n,$ and $H$ has property ${\bf Q},$  then $|H| < r^{t-1}$. This may be reduced to the case that $H$ is primitive as permutation group. If $F(H) = 1$, then $H$ has a faithful complex representation of degree $t-1 < n$, so we have $$|H|= [H:F(H)] < r^{t-2}$$ by the minimality of the pair $(G,n).$ If $F(H) \neq 1,$ then we are in the (extended) affine case, and by Corollary 4.2, we have  $$|H| \leq 24^{\frac{t-1}{3}} <   r^{t-1}.$$ 

\medskip
Next, we will reduce, using Lemma 8.2,  to the case $F(G) = Z(G).$ 

\medskip
For suppose otherwise. Then we may choose a nilpotent normal subgroup $M$ of $G$ such that $M/Z(G)$ is a minimal normal subgroup  of $G/Z(G)$. Then $M$ is non-Abelian,  and $M = Z(G)Q$ for some extra-special $q$-group $Q$ of order $q^{2m+1}$, where $Q \lhd G$, and $q^{m}$ is a divisor of $n$.

\medskip
Then ${\rm Out}_{G}(M)$ is isomorphic to a subgroup (with irreducible action on $Q/Z(Q))$ of the symplectic group  ${\rm Sp}(2m,q).$

\medskip
We have $$|{\rm Out}_{G}(M)|\leq |{\rm Sp}(2m,q)|  \leq {\rm s}(q,m)^{q^{m}-1} \leq 720^{\frac{q^{m}-1}{3}},$$ using Corollary 4.2.
Since $r > 720^{\frac{1}{3}},$ we certainly have $$|{\rm Out}_{G}(M)| <  r^{q^{m}-1}.$$

\medskip
Let us write $n =q^{m}f$ for an integer $f$. From Lemma 8.2, we know (since $M = F(M)$) that 
$$[G:F(G)] \leq |{\rm Out}_{G}(M)|[C_{G}(M): F(C_{G}(M))] .$$ If $C_{G}(M)$ is nilpotent, then we obtain 
$[G:F(G)] < r^{n-1},$ contrary to hypothesis. Thus $C_{G}(M)$ is certainly non-Abelian, so $f > 1$ and $n \geq 4.$

\medskip
Since $C_{G}(M)$ has property ${\bf Q},$  and $C_{G}(M)$ is isomorphic to a subgroup 
of ${\rm GL}(f,\mathbb{C})$ where $1 < f < n,$ the minimality of the pair $(G,n)$ tells us that $[C_{G}(M):F(C_{G}(M))] < r^{f-1}.$

\medskip
For any integer $u >4$ and any integer $1 < x < u$ which divides $u,$ we have $x + \frac{u}{x} \leq 2+\frac{u}{2},$ since the function
$x + \frac{u}{x}$ decreases monotonically on $(2,\sqrt{u})$ and increases monotonically on $(\sqrt{u},\frac{u}{2}).$

\medskip
Hence we have $$[G:F(G)] <  r^{q^{m}+f -2} \leq r^{\frac{n}{2}} \leq r^{n-1},$$ contrary to the choice of the pair $(G,n).$

\medskip
This contradiction shows that we have $F(G) \leq Z(G).$ In particular, the group $G$ in our minimal pair $(G,n)$ is not solvable.

\medskip
We now know that the group $G$ is primitive and that $F(G) = Z(G)$. There is a normal subgroup $M$ of $G$ with $Z(G) \leq M$ such that  $M/Z(G)$ is a minimal normal subgroup of $G/Z(G).$ By the minimal choice of $(G,n)$ we now have $Z(G) \leq \Phi(G) \leq F(G) = Z(G),$ so that 
$Z(G) = \Phi(G)$ (otherwise, we could write $G = Z(G)H$ for a maximal subgroup $H$, and then $H \lhd G$. In that case, $H$ 
also has property ${\bf Q}$ and $$[H:F(H)] = [H:Z(H)] = [G:Z(G)] \geq r^{n-1}).$$

\medskip
We next claim that $G$ has a unique component.

\medskip
By properties of the generalized Fitting subgroup $F^{\ast}(G),$ it is clear that $M = Z(G)E$, where $E$ is a central product of $t$ components, all $G$-conjugate, and all with the same centre $Z(E)$, say these are $\{L_{i}: 1 \leq i \leq t \}.$ 

\medskip
We next prove that $t = 1$ and $n = d.$ In particular, the group $E$ is quasi-simple and acts irreducibly.

\medskip
For let $$K = \bigcap_{i=1}^{t} N_{G}(L_{i}) \lhd G.$$ Then $G/K$ has property ${\bf Q},$  and is isomorphic to a subgroup of the symmetric group ${\rm Sym}_{t}$. Let $L_{i}/Z(E) \cong X$ (the same $X$ for each $i$), a finite non-Abelian simple group in  ${\rm s}({\bf Q}).$ By the structure of the automorphism groups of the finite simple groups, $K/EC_{G}(E)$ is solvable. Furthermore, $$[G:MC_{G}(M) \leq |{\rm Out}(X) \wr (G/K)|.$$

\medskip
By the primitivity of $G$, there is an integer $d \geq 2$ such that $E$ is isomorphic to an irreducible subgroup of $${\rm GL}(d^{t},\mathbb{C}),$$
each $L_{i}$ is isomorphic to an irreducible subgroup of ${\rm GL}(d,\mathbb{C}),$ and $d^{t}$ is a divisor of $n$. 
Certainly $$n \geq d^{t}  \geq 1+t(d-1) > t,$$ using the binomial theorem. Thus, by the minimality of the pair $(G,n)$, we have 
 $[G:K] <  r^{t-1}$, since $G/K$ is isomorphic to a subgroup of ${\rm Sym}_{t}$ ( and arguing as we did in the reduction to the primitive case earlier).

\medskip
We first prove that $E$ is irreducible. Otherwise, we have $n \geq 2d^{t}.$ In that case, we have $$[L_{i}:Z(E)] = |X|$$ for each $i$. We know that $X$ is generated by two elements, so that $$|{\rm Aut}(X)| < |X|^{2}$$ and $|{\rm Out}(X)| < |X|.$ Since $L_{i}$ is isomorphic to a subgroup of 
${\rm GL}(d, \mathbb{C})$ and $n \geq 2d^{t}$, we have $|X| = [L_{i}:Z(E)] <  r^{d-1}$ by the minimal choice of $(G,n).$

\medskip
Now we have $$[E:Z(E)] < r^{dt-t}$$ and also $$[K:C_{G}(E)] \leq |{\rm Aut}(X)|^{t} < |X|^{2t} \leq r^{2dt-2t}.$$  Hence 
$$[G:C_{G}(E)] < r^{2dt-t-1}.$$ 

\medskip
By Lemma 8.2, $C_{G}(E)$ is isomorphic to a subgroup of ${\rm GL}(\frac{n}{d^{t}}, \mathbb{C}),$ so that $$[C_{G}(E):F(C_{G}(E))] < r^{\frac{n}{d^{t}}-1},$$ by the minimality of $(G,n).$
But $$F(C_{G}(E)) \leq F(G) = Z(G) \leq F(C_{G}(E)),$$ so that now $$[G:Z(G)] <  r^{\frac{n}{d^{t}}  + 2dt - t - 2}.$$

\medskip
We have $$d^{t-1} \geq 1 + (t-1)(d-1) \geq t,$$ so that $d^{t} \geq dt$.

\medskip
Now we have $$\frac{n}{d^{t}} + 2dt-t-2 \leq \left(\frac{n}{d^{t}} + d^{t}\right) + dt-t -2 \leq \frac{n}{2} + dt-t \leq n-t \leq n-1,$$ contrary to the choice of $(G,n)$. Hence $E$ is irreducible, as claimed, so that $n = d^{t}.$

\medskip
If now $t > 1,$ then we may argue as above that $$[E:Z(G)] < r^{dt-t}$$ since $d < n$, and also that $$[K:Z(G)] < [E:Z(G)]^{2}.$$
Hence $$r^{d^{t}-1 } < [G:Z(G)] \leq r^{2dt-2t + t-1}.$$ Thus we have $d^{t}  <  2dt-t$. If $t = 2$, then $d^{2} < 4d$, so that $d \leq 3.$ 
If $t > 2,$ we have $$d^{t} \geq 1 + t(d-1) + t(d-1)^{2} = 1+td(d-1),$$
so that $$td(d-1)	\leq d^{t}-1 < 2dt-t.$$ Then $td(d-1) \leq  2t(d-1)$ and $d = 2.$

\medskip
Now $$2^{t} = d^{t} < 2dt-t = 3t,$$ which forces $t \leq 3.$ Hence only three possibilities remain:

\medskip
\noindent a) $d= 2, t = 3,$ and $n = 8$. For the sake of convenience, we invoke the fact that the finite subgroups of ${\rm GL}(2,\mathbb{C})$ are known. In this case, $E$ is the central product of three copies of ${\rm SL}(2,5)$, and ${\rm SL}(2,5) \in {\rm as}({\bf Q})$. By the minimality of $(G,n)$, we have  $60 \leq r^{2-1},$ so $ r \geq  60.$ Now $$[G:E] \leq |(\mathbb{Z}/2\mathbb{Z}) \wr {\rm Sym}_{3}|,$$ and so
$$[G:Z(G)] \leq  48 \times 60^{3} < 60^{4} < 60^{n-1} = 60^{7},$$ contrary to the choice of $(G,n).$

\medskip
\noindent b) $t = 2, d = 2,$ and $n = 4$. In this case, $E$ is the central product of two copies of ${\rm SL}(2,5)$. As above, $ r > 60.$ Now 
$$[G:E] \leq |(\mathbb{Z}/2\mathbb{Z}) \wr {\rm Sym}_{2}|,$$ and 
$$[G:Z(G)] \leq  8 \times 60^{2} < 60^{3} = 60^{n-1},$$ contrary to the choice of $(G,n).$

\medskip
\noindent c) $t = 2, d = 3,$ and $n = 9$. Again for the sake of convenience, we invoke the fact that the finite subgroups of ${\rm GL}(3,\mathbb{C})$ are long known.

\medskip
This time, $E = L_{1} L_{2}$, where $L_{1},L_{2}$ are distinct conjugate components of $G$, each isomorphic to one of ${\rm Alt}_{5}, {\rm PSL}(2,7)$ or a $3$-fold central extension  of ${\rm Alt}_{6}.$

\medskip
In each case, we obtain $$[G:Z(G)] \leq  32 \times [E:Z(E)], $$  and also $$r \geq \sqrt{|X|}$$ 
since $L_{1}/Z(L_{1}) \in s({\bf Q})$ and $L_{1}$ is isomorphic to a subgroup of ${\rm GL}(3,\mathbb{C}).$ But now 
$$r^{n-1} = r^{8} \geq  |X|^{4} > [G:Z(G)],$$ contrary to the choice of $(G,n).$

\medskip
This contradiction shows that we have $t = 1$, so that $G/Z(G)$ is almost simple, and $E$ is quasisimple and irreducible with $E/Z(E) \cong X.$
Also, $G/EZ(G)$ is now isomorphic to a subgroup of ${\rm Out}_{G}(X)$ and, in particular,\\ $[G:Z(G)] <|X|^{2}.$ Thus we must have 
$|X| > r^{\frac{n-1}{2}}$, since $[G:Z(G)] \geq r^{n-1}.$

\medskip
Suppose now that $E$ is imprimitive. It is convenient to work with characters for a while. Let $\chi$ be the character afforded by the given representation of $E$. Then $\chi = {\rm Ind}_{T}^{E}(\mu)$ for some irreducible character $\mu$ of a proper subgroup $T$ of $E$. 

\medskip
Then $n = [E:T]\mu(1)$ and $X$ has a faithful complex irreducible character of degree at most $[E:T]-1.$  
If $h \leq n-1$ is the minimum degree of a non-trivial complex irreducible character of $X$, then (as $X$ is simple), we have $|X|< r^{h-1}$
and then $$r^{n-1} \leq [G:Z(G)] < r^{2h-2}.$$ Thus  $2h-2 > n-1$, so that $h > \frac{n+1}{2}.$ In particular, $\mu$ must be linear.

\medskip
Furthermore, $X$ must act doubly transitively on the cosets of $T/Z$, for otherwise $X$ has a faithful complex irreducible character of degree at most $\frac{n-1}{2},$ contrary to the discussion above. 

\medskip
More precisely, we note that $$r^{n-1} \leq  [G:Z(G)] = |X||{\rm  Out}_{G}(X)|$$ and we know that $$[X:F(X)] = |X| < r^{n-2}.$$ 
For  $X \in {\rm s}({\bf Q})$ and $X$ is isomorphic to a subgroup of ${\rm GL}(n-1,\mathbb{C})$, so we may invoke the minimality of the pair
$(G,n).$ Then we have $|{\rm Out}_{G}(X)| > r >  8.$ No finite simple alternating group has an outer automorphism group of order greater than 
$4$, and no sporadic finite simple group has an outer automorphism group of order greater than $2,$ so that $X$ is neither a sporadic simple group nor an alternating group.

\medskip
Using Corollary 1.4 of (Maroti [9]), we may derive a contradiction to the assumed imprimitivity of $E$ . For $X$ is doubly transitive, hence certainly primitive as permutation group. If $|X| \leq 2^{n-1}$, then we have $[G:Z(G)] < 4^{n-1} < r^{n-1},$ a contradiction. If $|X| > 2^{n-1}$, then (since $X$ is not an alternating group) the possibilities for $X$ are restricted to the simple groups appearing in Corollary 1.4 of [9]. However, only three of those simple groups are neither alternating nor sporadic (recalling that ${\rm GL}(4,2)$ is isomorphic to ${\rm Alt}_{8}$), and these are ${\rm PSL}(2,7), {\rm PSL}(3,3)$ and ${\rm SL}(2,8)$.  These groups have outer automorphism groups of order $2,2$ and $3$ respectively, whereas $|{\rm Out}_{G}(X)| \geq r > 8,$ a contradiction. Hence the action of $E$ (as linear group) is primitive.

\medskip
By definition of the least upper bound $\ell({\bf Q}),$ we now have $$r^{n-1}  \leq [G:Z(G)] \leq \ell({\bf Q})^{n-1}$$
and $r \leq  \ell({\bf Q}),$ a contradiction.

\medskip
We conclude that there is no such pair $(G,n)$ with $r > {\rm max}\{720^{\frac{1}{3}}), \ell({\bf Q}) \},$ 
and that we have 
$$[G:F(G)] \leq  \left( {\rm max} \{ 720^{\frac{1}{3}},\ell({\bf Q}) \} \right)  ^{n-1}$$ whenever $G$ is a finite subgroup of
${\rm GL}(n,\mathbb{C})$ with property ${\bf Q}.$ In other words, we have 
$\alpha({\bf Q}) \leq {\rm max} \{720^{\frac{1}{3}}, \ell({\bf Q}) \}.$

\medskip
\noindent {\bf Remark 10.2 :}  Strictly speaking, we should repeat the arguments with $24$ in place of $720^{\frac{1}{3}}$ to obtain  
$\gamma({\bf Q}) \leq {\rm max} \{24, \ell({\bf Q}) \}.$ Also, we may argue as we did during the proof to show that, once $\alpha({\bf Q})$ is shown to exist, we have  $$\beta({\bf Q}) \leq {\rm max}\{720^{\frac{1}{3}},\alpha({\bf Q}) \}.$$ 

\medskip
\noindent {\bf Acknowledgements:} We are grateful to P. Etingof for the questions which initiated this work. We also thank Gunter Malle for a careful reading of an earlier version of this note, and for comments which improved the exposition, as well as correcting typographical errors and other inaccuracies, and Radha Kessar for further corrections. We also thank an anonymous referee for a thorough, constructive and helpful report which prompted a major revision of an earlier version of this note.

\section{References}

\medskip
\noindent [1] Blichfeldt, H.,  \emph{Finite Collineation Groups},University of Chicago Press, Chicago, (1917).

\medskip
\noindent [2] Collins, Michael J.,  \emph{On Jordan's theorem for complex linear groups}, Journal of  Group Theory, 10, (2007), 411-423.

\medskip
\noindent [3] Collins, Michael J.,  \emph{Bounds for primitive complex linear groups}, Journal of Algebra, 319, (2008), 759-776.

\medskip
\noindent [4] Coulembier, K.,  Etingof, P. and Ostrik, V., \emph{Asymptotic Properties of Tensor Powers in Symmetric Tensor Categories}, arXiv:2301.09804, January 2023.

\medskip
\noindent [5] Conway, J.H., Curtis, R.T., Norton, S.P.,  Parker, R.A and  Wilson, R.A., \emph{Atlas of Finite Groups}, Clarendon, Oxford, (1985).

\medskip
\noindent [6] Dornhoff, L., \emph{Group Representation Theory, Part A : Ordinary Representation Theory}, Marcel Dekker, New York, (1972).

\medskip
\noindent [7] Jansen, C., \emph{The minimal degrees of the faithful representations of the sporadic simple groups and their covering groups}, 
LMS Journal of Computational Mathematics,8,(2005),122-144.

\medskip
\noindent [8] Landazuri, V., Seitz, G.M., \emph{On the minimal degree of projective representations of finite Chevalley groups}, Journal of Algebra, 32,2,(1974), 418-443.

\medskip
\noindent [9] Maroti, A., \emph{On the order of primitive groups}, Journal of Algebra, 258, (2002), 631-640.

\medskip
\noindent [10] Robinson, G.R., \emph {Bounding the Size of Permutation Groups and Complex Linear Groups of Odd Order}, Journal of Algebra, 335, (2011), 163-170.

\medskip
\noindent [11] Tiep, P.H., \emph{Low dimensional representations of finite quasisimple groups} {\bf in} \emph{Groups, Combinatorics and Geometry, Durham 2001}, World Scientific, River Edge, New Jersey, 2003.

\medskip
\noindent [12] Weisfeiler, B., \emph{Post-classification version of Jordan's Theorem on finite linear groups}, Proceedings of National Academy of Science, 81, 16, (1984), 5278-5279.

\medskip
\noindent [13] Weisfeiler, B., \emph{On the size and structure of finite linear groups}, arXiv 1203.1960v2 .

\end{document}